\documentclass[12pt]{article} 
\usepackage{epic,eepic}
\usepackage{epsfig}
\newtheorem{theor}{Theorem}
\begin{document}

\title{
Embedding the graphs of regular tilings and star-honeycombs into the graphs
of hypercubes and cubic lattices
\thanks{This work was supported by the Volkswagen-Stiftung 
(RiP-program at Oberwolfach) and Russian fund of fundamental 
research (grant 96-01-00166).}}

\author{Michel DEZA \\
     CNRS and Ecole Normale Sup\'erieure, Paris, France \\
\and   Mikhail SHTOGRIN \\
Steklov Mathematical Institute, 117966 Moscow GSP-1, Russia} 

\date{}

\maketitle 

\begin{abstract}
We review the regular tilings of $d$-sphere, Euclidean $d$-space, hyperbolic
$d$-space and Coxeter's regular hyperbolic honeycombs (with infinite or
star-shaped
cells or vertex figures) with respect of possible embedding, isometric up to a
scale, of their skeletons into a m-cube or m-dimensional cubic lattice.
In section 2 the last remaining 2-dimensional case is decided: for any odd
$m \ge 7$, star-honeycombs $\{m, \frac {m}{2}\}$ are embeddable while
$\{ \frac {m}{2}, m\}$ are not (unique case of non-embedding for dimension 2).
As a spherical analogue of those honeycombs, we enumerate, in section 3, 36
Riemann surfaces representing all nine regular polyhedra on the sphere. In
section 4, non-embeddability of all remaining star-honeycombs (on
3-sphere and hyperbolic 4-space) is proved. In the last section 5, all cases
of embedding for dimension $d > 2$ are identified. Besides hyper-simplices and
hyper-octahedra, they are exactly those with bipartite skeleton:
hyper-cubes, cubic lattices and 8, 2, 1 tilings of hyperbolic 3-, 4-,
5-space (only two, $\{4,3,5\}$ and $\{4,3,3,5\}$, of those 11 have compact both, facets
and vertex figures).

\end{abstract}
\section{Introduction}

We say that given tiling (or honeycomb) $T$ has a {\em $l_1$-graph}
and embeds up to {\em scale} $\lambda$ into $m$-cube $H_m$ (or, if the
graph is infinite, 
into  cubic lattice ${\bf Z}_m$ ), if there exists a mapping $f$ of 
the vertex-set of the skeleton graph of $T$ into the vertex-set of
$H_m$ (or ${\bf Z}_m$) such that
\[\lambda d_{T}(v_i,v_j)=||f(v_i),f(v_j)||_{l_1} =
\sum_{1 \le k \le m}|f_{k}(v_i)-f_k(v_j)| \mbox{ for all vertices
$v_i,v_j$}, \]
where $d_{T}$ denotes the graph-theoretical distance in contrast to the
normed-space distance $l_1$.
The smallest such number $\lambda$ is called {\em minimal scale}.

Denote by $T \to H_m$ (by $T \to {\bf Z_m}$) isometric
embedding of the skeleton graph of $T$ into $m$-cube (respectively, into
$m$-dimensional
cubic lattice); denote by $T \to \frac{1}{2} H_m $ and by 
$T \to \frac{1}{2} {\bf Z_m}$ isometric up to scale 2 embedding.

Call an embeddable tiling $l_1$-{\em rigid}, if all its embeddings as above 
are pairwise equivalent.
All, except hyper-simplexes and hyper-octahedra (see Remark 4  below), embeddable tilings in this paper turn out to be $l_1$-rigid
and so, by a result from \cite{Shp}, having scale 1 or (only for non-bipartite planar tilings) 2. Those
embeddings were obtained by constructing a complete system of
{\em alternated zones }; see \cite{CDG}, \cite{DS1}, \cite{DS2}.   

Actually, a tiling is a special case of a honeycomb, but we reserve the
last term for the case  when the cell or the vertex figure is a
star-polytope and so the honeycomb covers the space several times; the
multiplicity of the covering is called its {\em density}. 
 
Embedding of Platonic solids was remarked in \cite{Ke} and precised,
for the dodecahedron, in \cite{AD}. Then \cite{As} showed that
$\{3,6\}, \{6,3\}$,
and $\{m,k\}$ (for even $m \ge 8$ and $m = \infty $) are embeddable. The remaining
case of odd $m$ and limit cases of $m = 2, \infty $ was decided in \cite{DS1};
 all those results are put together in the Theorem 1 below.

All four star-polyhedra are embeddable.
The great icosahedron $\{3,\frac{5}{2}\}$ of Poinsot and the great stellated
dodecahedron  $\{ \frac{5}{2}, 3\}$ of Kepler have the skeleton (and, moreover,
the surface) of, respectively, icosahedron and dodecahedron; each of them
has density 7. All ten star-4-polytopes are not embeddable: see the
Theorem 3 below. 

The case of Archimedean tilings of 2-sphere and of Euclidean plane was
decided in \cite{DS1}; it turns out that for any such tilings (except
{\em $Prism_3$} and its dual, both embeddable) exactly one of two (a tiling
and its dual) is embeddable. For 3-sphere and 3-space it was done in
\cite{DS4}; for example, Gosset's semiregular 4-polytope snub 24-cell turns
out to be embeddable into half-12-cube.
All 92 {\em regular-faced} 3-polytopes were considered in \cite{DG2} 
and, for all higher dimensions, in \cite{DS1}. The {\em truncations}
of regular polytopes were considered in \cite{DS2}. Another large
generalization of Platonic solids - {\em bifaced} polyhedra - were
considered in \cite{DG2}. (Some generalizations of Archimedean plane
tilings, {\em 2-uniform} ones and  {\em equi-transitive} ones, were treated
in \cite{DS1}, \cite{DS2}, respectively.) Finally, skeletons of
(Delaunay and Voronoi tilings of) lattices were dealt with 
in \cite{DS3}.
 
Embeddable ones, among all regular tilings of all dimensions, having compact
facets and vertex fugures, were identified in \cite{DS1}, \cite{DS2}.
 
Coxeter (see \cite{Cox}) extended the concept of regular tiling, permitting
infinite cells and vertex figures, but with the fundamental region of the 
symmetry group of a finite content. His second extension was to
permit {\em honeycombs}, i.e. star-polytopes can be cells or
vertex figures. For the 2-dimensional case, on which we are focusing in the 
next Section, above
extensions produced only following new honeycombs - $\{ \frac{m}{2}, m\}$ and
$\{m, \frac{m}{2}\}$ for any odd $m \ge 7 $ - which are hyperbolic analogue
of spherical star-polyhedra $\{ \frac{5}{2}, 5\}$ (the small stellated dodecahedron
of Kepler) and $\{5, \frac{5}{2}\}$ (the great dodecahedron of Poinsot). Both
$\{ \frac{5}{2}, 5\}$ and $\{5, \frac{5}{2}\}$ have the skeleton of the icosahedron. For
any odd $m$ above honeycombs cover the space (2-sphere for $m = 5$) 3 times.
The skeleton of $\{m, \frac{m}{2}\}$ is, evidently, the same as of 3m, because it 
can be seen as $\{3,m\}$ with the same vertices and edges forming $m$-gons instead of
triangles. The faces of $\{ \frac {m}{2}, m\}$ are stellated faces of $\{m,3\}$
and it have the same vertices as $\{3,m\}$.

We adopt here classical definition of the regularity: the transitivity of the
group of symmetry on all faces of each dimension. But, as remarked the
referee, the modern concept of regularity, which requires transitivity on
flags, would not necessitate any change in the results.

The following {\em 5-gonal} inequality (see \cite{D}):

$d_{ab} + (d_{xy} + d_{xz} + d_{yz}) \le (d_{ax} + d_{ay} + d_{az}) +
(d_{bx} + d_{by} + d_{bz})$

for distances between any five vertices $a, b, c, x, y$, is an important
 necessary condition for embedding of
graphs, which will be used in proofs of Theorems 3,4 below.

This paper is a continuation of general study of $l_1$-graphs
and $l_1$-metrics, surveyed in the book \cite{DL}, where many applications and
connections of this topic are given. In addition, we tried here to extract
from purely geometric, affine structures, considered below, their new, purely
combinatorial (in terms of metrics of their graphs) properties.

\section{ Planar tilings and hyperbolic honeycombs}
 
They are presented in the Table 1 below; we use the following notation:

1. The row indicates the facet (cell) of the tiling (or honeycomb),
the column indicates its vertex figure. The tilings and honeycombs are
denoted usually by their
Schl\"afli notation, but in the Tables 1, 3-5 below we omit the brackets and
commas for convenience (in order to fit into page). 

2. By $m$  we denote m-gon, by $\frac{m}{2}$ star-m-gon (if 
$m$ is odd). By
$ \alpha_3$, $ \beta_3$, $ \gamma_3$, Ico, Do and $ \delta_2$ we denote
regular ones tetrahedron $\{3,3\}$, octahedron $\{3,4\}$, cube $\{4,3\}$,
icosahedron $\{3,5\}$,
dodecahedron $\{5,3\}$ and the square lattice $Z_2 = \{4,4\}$. The numbers are: any
$ m\ge 7$ in 8th column, row and any {\em odd } $ m\ge 7$ in 9th column, row.

3. We consider that: $\{2,m\}$ is a 2-vertex multi-graph with $m$ edges; $\{m,2\}$
can be seen as a $m$-gon; all vertices of $m \infty $ are on the absolute
 conic at infinity (it has an infinite degree); the faces $ \infty $ of
 $ \{\infty, m\}$ are inscribed in {\em horocycles} (circles centered in $\infty$).

{\bf Table 1. 
 2-dimensional regular tilings and honeycombs.}
\[\begin{array}{|c||c|c|c|c|c|c|c|c||c||c||} \hline
 &2&3 &4 & 5 &6 &7&m& \infty &  
\frac{m}{2}&\frac{5}{2}\\ \hline \hline
2&22&23&24&25&26&27&2m&2 \infty&& \\ \hline  
3&32&\alpha_3&\beta_3&Ico&{\bf 36}&37&3m&3 \infty & &3 \frac{5}{2} \\
 \hline  
4&42&\gamma_3& {\bf \delta_2 } &45&46&47&4m&4 \infty && \\ \hline        
5&52&Do&54&55&56&57&5m&5 \infty & &5 \frac{5}{2}    \\ \hline
6&62&{\bf 63}&64&65&66&67&6m&6 \infty  && \\ \hline 
7&72&73&74&75&76&77&7m&7 \infty &&  \\ \hline
m&m2&m3&m4&m5&m6&m7&mm&m \infty&m \frac {m}{2}  & \\ \hline
 \infty & \infty 2& \infty 3& \infty 4& \infty 5& \infty 6& \infty 7&
 \infty m& \infty \infty && \\ \hline \hline
\frac{m}{2}&&&&&&& \frac{m}{2} m&&& \\ \hline \hline
\frac{5}{2}&&\frac{5}{2} 3&&\frac{5}{2} 5&&&&&& \\ \hline \hline
\end{array} \]

\begin{theor}
\label{2til}
 All 2-dimensional tilings $\{m,k\}$ are embeddable,namely:

(i) if $ \frac {1}{m} + \frac {1}{k}  >  \frac {1}{2} $ ( 2-sphere),
 then

$\{2,m\} \to H_1$ for any $m$, $\{m,2\} \to \frac{1}{2} H_m$ for odd $m$ and
$\{m,2\} \to H_{\frac{m}{2}}$ for even $m$;

 $ \{3,3\}=\alpha_3  \to \frac{1}{2} H_3$, $\frac{1}{2} H_4$; 
  $ \{4,3\}=\gamma_3  \to H_3$;  $ \{3,4\}=\beta_3  \to \frac{1}{2} H_4$;

$ \{3,5\}=Ico ( \sim \{3, \frac{5}{2}\} \sim \{5, \frac{5}{2}\} \sim \{\frac{5}{2}, 5\})
 \to H_6$  and  $ \{5,3\}=Do ( \sim \{\frac{5}{2}, 3\})  \to \frac{1}{2} H_{10}$;

(ii) if $ \frac {1}{m} + \frac {1}{k} = \frac {1}{2}$ (Euclidean plane), then

$\{2, \infty\} \to H_1$, $ \{\infty, 2\} \to Z_1$;  $\{4,4\}= \delta_2 \to {\bf Z}_2$, $\{3,6\} \to \frac{1}{2} {\bf Z}_3$, $\{6,3\} \to {\bf Z}_3$; 

(iii) if $ \frac {1}{2}  >  \frac {1}{m} + \frac {1}{k} $ (hyperbolic plane),
then

 $\{m,k\}  \to \frac{1}{2} {\bf Z}_ \infty $ if m is odd, $ k \le \infty$
 and $\{m,k\}  \to  {\bf Z}_ \infty $ is $m$ is even or $ \infty $, 
$k \le \infty$.

\end{theor}

{\bf Remark 1} (notation and terms here are from \cite{Cox3}, \cite{Cro}): 

(i) The embedding of the icosahedron  $\{3,5\}$ into $\frac{1}{2}H_6$ was mentioned
in \cite{Cox2} on pages 450--451, as {\em regular skew icosahedron}. There are
5 proper regular-faced fragments of $\{3,5\}$: 5-pyramid, 5-antiprism,
para-bidiminished $\{3,5\}$, meta-bidiminished $\{3,5\}$, and tridiminished
$\{3,5\}$; 
5-pyramid embeds into $\frac{1}{2}H_5$, all others into $\frac{1}{2}H_6$.

(ii) The antipodal quotients of (embeddable, see Theorem 1 (i) above)
cube, icosahedron, dodecahedron are regular maps
$\{4, 3\}_3, \{3, 5\}_5, \{5, 3\}_5$ on the projective plane, which are
$K_4, K_6$, the Petersen graph; they embed into $\frac{1}{2}H_m$
for $m = 4, 6, 6$, respectively.

(iii) Besides $\{4,4\}, \{3,6\}, \{6,3\}$ (embeddable, see Theorem 1 (ii) above), there are 
exactly three other infinite regular polyhedra. They are {\em regular skew}
polyhedra $\{4, 6|4\}$, $\{6, 4|4\}$, $\{6, 6|3\}$, which can be obtained by
deleting of cells from the tilings of 3-space by cubes (${\bf Z}_3$), by
truncated octahedra (the Voronoi tiling for the lattice $A^*_3$), by regular
tetrahedra and truncated tetrahedra (Föppl uniform tiling). They are,
respectively: embeddable into ${\bf Z}_3$, embeddable into ${\bf Z}_6$,
not 5-gonal. All finite regular skew 4-polytopes are: the family $\{4, 4|m\}$
of self-dual quadringulations of the torus (it is the product
of two $m$-cycles and so embeddable into $\frac{1}{2}H_{2m}$ for odd $m$ or
into $H_m$ for even $m$), not 5-gonal $\{6, 4|3\}, \{4, 6|3\}, \{8, 4|3\}$ and
its undecided dual $\{4, 8|3\}$. 

(iv) Examples of other interesting regular maps are the Dyck map $\{3,8\}_6$
(8-valent map with 12 vertices and 32 triangular faces),the Klein map
$\{3,7\}_8$ (7-valent map with 24 vertices and 56 triangular faces) and
$\{4,5\}_5$ (5-valent map with 16 vertices and 20 quadrangular faces). Those
maps (all of oriented genus 3) come from the hyperbolic tilings $\{3,8\},
\{3,7\}, \{5,4\}$,
respectively (which are embeddable; see Theorem 1 (iii) above) by 
identification of some vertices of the unit cell. Those three maps and their
duals are all not 5-gonal. But, for example, the 3-valent partition of the
torus into 4 hexagons is embeddable: it is the cube on the torus.

{\bf Remark 2} (notation and terms here are from \cite{Cox1}, \cite{Wen} and
\cite{Cro}). With V.P.Grishukhin we considered embeddability of following
non-convex polyhedra:

(i) All non-Platonic facetings of Platonic solids (see
\cite{Cox1}, page 100) are: 4 star-polyhedra $\{ \frac{5}{2},5\} $,
$\{5, \frac{5}{2}\}$,
$\{ \frac{5}{2},3\}$, $\{3, \frac{5}{2}\}$ and 4 regular compounds $2 \alpha_3$ (Kepler's
{\em stella octangula}), $5 \gamma_3$, $5 \alpha_3$, $10 \alpha_3$. The 
remaining regular compound is $5 \beta_3$, which is dual to $5 \gamma_3$. In
this Remark only, contrary to Theorem 1 (i), we consider all visible 
``vertices'' of polyhedra, not only those of their shells. Then it turns out,
that $\{\frac{5}{2},5\}$, $\{5, \frac{5}{2}\}$, $\{\frac{5}{2},3\}$,
$\{3, \frac{5}{2}\}$,
$2 \alpha_3$, $5 \beta_3$ have the same skeletons as dual truncated,
respectively, $\{3,5\}$, $\{5,3\}$, $\{5,3\}$, truncated $\{3,5\}$, $\gamma_3$, icosidodecahedron.
$5 \alpha_3$ has the same skeleton as dual snub dodecahedron. Among all 4
star-polyhedra, 5 regular compounds and their 9 duals, all embeddable ones are:

$\{\frac{5}{2},5\} \to \frac{1}{2} H_{10}$,
$\{5, \frac{5}{2}\} ( \sim \{ \frac{5}{2},3\}) \to \frac{1}{2} H_{26}$,
$\{3, \frac{5}{2}\} \to \frac{1}{2} H_{70}$, $2 \alpha_3 \to \frac{1}{2} H_{12}$,
dual $5 \beta_3 ( \sim 5 \gamma_3) \to H_{15}$,
dual $5\alpha_3 \to \frac{1}{2}H_{15}$.

(ii) Among 8 stellations $A-H$ of $\{3,5\}$ (the main sequence, see \cite{Cro},
page 272), all embeddable ones are $A = \{3,5\}$, $B \sim \{5, \frac{5}{2}\}$
and $G \sim H \sim \{3, \frac{5}{2}\}$. Also the dual of the stellation $De_2f_2$ of
$\{3,5\}$ has the same skeleton as the rhombicosidodecahedron and it embeds into
$ \frac{1}{2} H_{16}$. The stellations $De_1 \sim Fg_2 \sim C =5 \alpha_3$ and
$Fg_1$, $De_2f_2$ are not embeddable.

iii) Among the compounds of two dual Platonic solids and dual compounds, all
embeddable ones are $2 \alpha_3$ and, into $ \frac{1}{2} H_{28}$, the dual of
$\{3,5\} + \{5,3\}$. Among all 53 non-convex non-regular uniform polyhedra (Nos.
67--119 in \cite{Wen}), two are
quasi-regular: the dodecadodecahedron and the great icosi
dodecahedron
(see \cite{Cox1}, page 101 and Nos. 73, 94 in \cite{Wen}). Again we consider
all visible ``vertices'' and see a pentagram $\frac{5}{2}$ as {\em pentacle}
(10-sided non-convex polygon). Then both above polyhedra and their duals are
not embeddable. But, for example, the ditrigonal dodecahedron (No. 80 in
\cite{Wen}, a relative of No. 73) embeds into $\frac{1}{2}H_{20}$.

The following theorem gives the family of all non-embeddable regular planar
cases.

\begin{theor}
\label{2hon}
For any odd $m \ge 7$ we have

(i) $\{ \frac{m}{2}, m\}$ is not embeddable;

(ii) $\{m, \frac{m}{2}\} ( \sim \{3,m\} ) \to \frac{1}{2} {\bf Z}_ \infty $.

\end{theor}

The assertion (ii) is trivial. The proof of (i) will be
preceded by 3 lemmas and first two of them are easy but of independent
interest for embedding of (not necessary planar) graphs. Lemma 1 can be
extended on the isometric cycles. 

Let $G$ be a graph, scale $ \lambda $ embeddable into ${ \bf Z_m}$, let $C$
be an oriented circuit of length $t$ in $G$ and let $e$ be an arc in $C$.
Then there are $ \lambda $ {\em elementary vectors}, i.e. steps in the
cubic lattice ${ \bf Z_m}$, corresponding to the arc $e$; denote them by
$x_1(e), ...,x_{ \lambda (e)}$. Clearly, the sum of all vectors $x_i(e)$ by
all $i$ and all arcs $e$ of the circuit, is the zero-vector.

Now, if $t$ is even, denote by $e^*$ the arc opposite to $e$ in the circuit
$C$; if $t$ is odd, denote by $e'$,$e''$ two arcs of $C$ opposite to $e$.
For even $t$, call the arc $e$ { \em balanced} if the  set of all its vectors
$x_i(e)$ coincides with the set of all $x_i(e^*)$,
but the vectors of arc $e^*$ go in opposite direction
{ \em on the circuit $C$ } to the vectors of $e$. For odd $t$, call the arc
$e$ { \em balanced} if a half of vectors of $e'$ together with a half of
 vectors of the second opposite arc $e''$ form a partition of the
set of vectors of $e$ and those vectors go in opposite direction (on $C$) to
those of arc $e$.  

Remind, that the {\em girth} of the graph is the length of its minimal
circuit.

\noindent
{\bf Lemma 1.}\ {\em Let $G$ be an embeddable graph of girth $t$. Then

(i) any arc of a circuit of length $t$ is balanced;

(ii) if $t$ is even, then any arc of a circuit of length $t + 1$ is also
balanced.}

\noindent
{\bf Lemma 2.}\ {\em Let $G$ be an embeddable graph of girth $t$ and let $P$ be
an isometric oriented path of length at most $ \lfloor \frac{t}{2} \rfloor $
in $G$. Then there are no two arcs on this path having vectors, which are
equal, but have opposite directions on the path.}

\noindent
{\bf Lemma 3.}\ {\em The girth of the skeleton of $ \{\frac{m}{2}, m\}$ is 3 for
$m = 5$ and $m - 1$ for any odd $m \ge 7$.}

\medskip

\noindent
{\bf Proof of Lemma 3}\

\medskip
\begin{center}
\setlength{\unitlength}{0.00066667in}
\begingroup\makeatletter\ifx\SetFigFont\undefined%
\gdef\SetFigFont#1#2#3#4#5{%
  \reset@font\fontsize{#1}{#2pt}%
  \fontfamily{#3}\fontseries{#4}\fontshape{#5}%
  \selectfont}%
\fi\endgroup%
{\renewcommand{\dashlinestretch}{30}
\begin{picture}(6109,4029)(0,-10)
\dashline{60.000}(3225,3054)(5325,3729)(5325,1404)
	(3225,1779)(4125,3804)(5925,2529)
	(4050,1029)(3225,3054)
\put(3150,1254){\makebox(0,0)[lb]{\smash{{{\SetFigFont{10}{12.0}{\rmdefault}{\mddefault}{\updefault}b3}}}}}
\put(3150,1554){\makebox(0,0)[lb]{\smash{{{\SetFigFont{10}{12.0}{\rmdefault}{\mddefault}{\updefault}a0}}}}}
\put(750,1029){\makebox(0,0)[lb]{\smash{{{\SetFigFont{10}{12.0}{\rmdefault}{\mddefault}{\updefault}a1}}}}}
\put(675,3804){\makebox(0,0)[lb]{\smash{{{\SetFigFont{10}{12.0}{\rmdefault}{\mddefault}{\updefault}a2}}}}}
\put(3150,3204){\makebox(0,0)[lb]{\smash{{{\SetFigFont{10}{12.0}{\rmdefault}{\mddefault}{\updefault}a3}}}}}
\put(3150,3504){\makebox(0,0)[lb]{\smash{{{\SetFigFont{10}{12.0}{\rmdefault}{\mddefault}{\updefault}b0}}}}}
\path(3225,3054)(825,3729)(825,1329)
	(3225,1779)(2175,3804)(150,2529)
	(2325,1029)(3225,3054)
\put(2175,804){\makebox(0,0)[lb]{\smash{{{\SetFigFont{10}{12.0}{\rmdefault}{\mddefault}{\updefault}a4}}}}}
\put(1425,54){\makebox(0,0)[lb]{\smash{{{\SetFigFont{10}{12.0}{\rmdefault}{\mddefault}{\updefault}           Fig . 1a.      A fragment of    7/2  7}}}}}
\put(0,2604){\makebox(0,0)[lb]{\smash{{{\SetFigFont{10}{12.0}{\rmdefault}{\mddefault}{\updefault}a5}}}}}
\put(1950,3879){\makebox(0,0)[lb]{\smash{{{\SetFigFont{10}{12.0}{\rmdefault}{\mddefault}{\updefault}a6}}}}}
\put(5400,3804){\makebox(0,0)[lb]{\smash{{{\SetFigFont{10}{12.0}{\rmdefault}{\mddefault}{\updefault}b1}}}}}
\put(5400,1254){\makebox(0,0)[lb]{\smash{{{\SetFigFont{10}{12.0}{\rmdefault}{\mddefault}{\updefault}b2}}}}}
\put(4125,3879){\makebox(0,0)[lb]{\smash{{{\SetFigFont{10}{12.0}{\rmdefault}{\mddefault}{\updefault}b4}}}}}
\put(5925,2679){\makebox(0,0)[lb]{\smash{{{\SetFigFont{10}{12.0}{\rmdefault}{\mddefault}{\updefault}b5}}}}}
\put(4125,879){\makebox(0,0)[lb]{\smash{{{\SetFigFont{10}{12.0}{\rmdefault}{\mddefault}{\updefault}b6}}}}}
\end{picture}
}

\end{center}
\medskip

Consider Fig. 1a. Take a cell $A = (a_0,...,a_m = a_0)$ of the
$ \{\frac{m}{2}, m\}$, i.e. a star 
$m$-gon, seen as an oriented cycle of length $m = 2k + 1$. Consider following
automorphism of the honeycomb: a turn by 180 degrees around the mid-point of
the segment $[ a_0, a_k]$. The image of $A$ is the oriented
star $m$-gon $B = (b_0,..., b_m = b_0)$ with $b_0 = a_k$, $b_k = a_0$.
Consider now oriented cycle $C = (a_0, a_1,...,a_k = b_0,..., b_k = a_0)$
of even length $m-1 = 2k$. In order to prove the Lemma 3, we will show that
$C$ is a cycle of minimal length.

First we show that the graph distance $d(a_0,a_k) = k$, i.e. the path
$P: = (a_0, a_1,...,a_k)$ is a shortest path from $a_0$ to $a_k$. It will
imply that $d(a_0, c(A)) = d(a_k, c(A)) = k$, where $c(A)$ is the center of
the cell $A$, because all vertices of $\{\frac{m}{2},m\}$ are vertices of regular
triangles of $\{3,m\}$.

Let $Q$ be a shortest path from $a_0$ to $a_k$. Then it goes around the vertex 
$c(A)$ or the center $c(B)$ of the cell $B$, because otherwise
$Q$ goes through at least one of the vertices $a_{k+1}$, $a_{2k}$, $b_{k+1}$,
$b_{2k}$ and so $Q$ contains at least one of the pairs of vertices
$(a_0, a_{k+1})$, $(a_0=b_k, b_{2k})$, $(b_k=a_0, a_{2k})$,
$(a_k=b_0, b_{k+1})$. But each of those pairs has, by the symmetry of our
honeycomb $\{ \frac{m}{2},m\}$, same distance between them as $(a_0, a_k)$; it
contradicts to the supposition that $Q$ is a shortest path. So, we can
suppose that $Q$ goes around $c(A)$ (the argument is the same if it goes
around $c(B)$). Now, to each edge $(s, t)$, corresponds, from
the center $c(A)$ of $A$, the angle $(s, c(A), t)$. The $2k+1$
edges of $A$ are only edges, for which this angle is $ \frac{4k \pi}{2k+1}$;
for any other edge, the angle is smaller, since it is more far from $c(A)$.
So, if $Q$ contains an edge, other than one from $A$, then, in order to reach
$a_k$ from $a_0$, it should be of length more than $k$. Therefore, any
shortest path from $a_0$ to $a_k$, should consist only of edges of $A$ and
then it is of length $k$. So, $d(a_0, c(A)) = k$ also, as well as for any edge of $\{3,m\}$. Same holds for $m = 5$.

We will show now that:

(i) any path $R$ of length $2k - 2$ is not closed and

(ii) $R$ cannot be closed by only one edge.

 But $C$ is a closed path of length
$2k$; so (i), (ii) will imply that $2k$ (respectively, $2k + 1$) is the
minimal length of any (respectively, any odd) simple isometric cycle in the
graph. For $m = 5$ (ii) does not holds.

Suppose that $R$ is closed; let as see it as a $2k - 2$-gon on hyperbolic
plane. Any angle of $R$ is a multiple $i \frac{2 \pi}{m}$, but $i > 1$
for at least one angle, because $(2k - 2) \frac{2 \pi}{m} < 2 \pi$.
Suppose that a angle has $1 < i \le k $; the argument will be the same if 
$k+1 \le i < m-1$, but
for the complementary angle $(m -i) \frac{2 \pi}{m} $ with $1 < m-i \le k$.

See Fig. 1b for the following argument. Fix an angle $r, s, t$ between two
adjacent edges $(r, s)$ and $(s, t)$ of
$R$. Let $s*$ be the opposite vertex to $s$ on $R$,
let $(s, r')$, $(s, t')$ be the edges such that the angles $r, s, r'$,
$t, s, t'$ are $\frac{2 \pi}{m}$. Let $A$, $B$ be the cells $\frac{m}{2}$,
defined by pairs $(r, s), (s, r')$ and $(t, s), (s, t')$ of their adjacent
edges and $c(A)$, $c(B)$ are their centers. The vertex $c(A)$  not belongs to
the path from $s$ to $s*$ of length $k - 1$, since we proved above that
$d(s, c(A) = k)$; so this path should go around $c(A)$. Let $p$ be the vertex
of $A$, reachable from $s$ by $k - 1$ steps on $A$, starting by $r$; let $q$
be the vertex of $B$, reachable from $s$ by $k -1$ steps on $B$, starting
by $t$. By mirror on $(r, s)$ (respectively, $(s, t)$) we obtain the cells
$A'$, $B'$, their centers $c(A')$, $c(B')$ and vertices $p'$, $q'$, which are
reflections of $p$, $q$. Call { \em A-domain}, the part of the hyperbolic
plane, bounded by half-lines ($c(A), p, \infty$), ($c(A'), p' \infty$) and the
angle ($c(A), s, c(A')$); call { \em B-domain}, the part, bounded similarly for
$B$. Actually, $B$-domain is the reflection of $A$-domain by the bisectrisse
of the angle $(r, s, t)$.
  
We will show now that the vertex $s*$ should belong to both $A$- and
$B$-domains. But they do not have common points, besides $s$. This
contradiction will show
that our $R$, a closed path of length $2k - 2$, do not exists. Any edge of
the path $(s, t, ...,s*)$ of length $k - 1$ is seen from $c(A)$ under angle
at most $ \frac{4 \pi}{m}$ with equality if and only if this edge belongs to
$A$ (as, for example, the edge $(r, s)$). Summing up those angles along the
path $(s t,...,s*)$, we get less than $(k-1) \frac{4 \pi}{m}$, obtained for
the path of length $k - 1$ from $s$ to $p$, going along $A$. It implies that
$s*$ belongs to $A$-domain and also, by reflection, to $B$-domain.

\medskip
\begin{center}
\setlength{\unitlength}{0.00041667in}
\begingroup\makeatletter\ifx\SetFigFont\undefined%
\gdef\SetFigFont#1#2#3#4#5{%
  \reset@font\fontsize{#1}{#2pt}%
  \fontfamily{#3}\fontseries{#4}\fontshape{#5}%
  \selectfont}%
\fi\endgroup%
{\renewcommand{\dashlinestretch}{30}
\begin{picture}(10074,9276)(0,-10)
\put(3853.071,5322.214){\arc{273.162}{3.9825}{7.4079}}
\put(5074.500,5461.500){\arc{530.330}{2.9997}{5.4978}}
\put(5862.000,6324.000){\arc{750.000}{2.2143}{3.1416}}
\dashline{60.000}(12,5049)(10062,5049)
\path(3312,5049)(1287,8949)
\path(3312,5049)(1287,1149)
\path(3312,5049)(9162,6399)
\path(3312,5049)(9237,3774)
\path(3312,5049)(3987,6549)
\path(3312,4974)(4062,3549)
\path(3312,5049)(4662,6099)(9012,7149)
\path(3312,5049)(4737,3924)(9087,3024)
\path(3312,5049)(3237,6624)(1887,9174)(1887,9249)
\path(3312,5049)(3237,3399)(1812,774)
\path(3987,6549)(5187,6999)(6087,6399)
\path(3987,6549)(3612,7524)(2712,7599)
\path(4062,3549)(5187,3099)(6162,3624)
\path(4062,3549)(3687,2499)(2712,2349)
\path(3312,5049)(5262,5499)
\put(4683.429,6088.286){\arc{441.762}{5.8847}{8.7425}}
\dashline{60.000}(5112,5499)(6012,6399)(6087,6399)
\put(3237,99){\makebox(0,0)[lb]{\smash{{{\SetFigFont{12}{14.4}{\rmdefault}{\mddefault}{\updefault}Fig. 1b.   A fragment of   9/2 9}}}}}
\dashline{60.000}(4662,6099)(5187,6999)
\dashline{60.000}(4662,6099)(3987,6549)
\dashline{60.000}(5037,4674)(4737,3849)(4812,3999)
\dashline{60.000}(3987,6549)(5037,8574)
\dashline{60.000}(4062,3474)(4062,3549)(5187,1449)
\put(2937,5124){\makebox(0,0)[lb]{\smash{{{\SetFigFont{6}{7.2}{\rmdefault}{\mddefault}{\updefault}s}}}}}
\put(3762,6474){\makebox(0,0)[lb]{\smash{{{\SetFigFont{6}{7.2}{\rmdefault}{\mddefault}{\updefault}r}}}}}
\put(3762,3474){\makebox(0,0)[lb]{\smash{{{\SetFigFont{6}{7.2}{\rmdefault}{\mddefault}{\updefault}t}}}}}
\put(6012,6549){\makebox(0,0)[lb]{\smash{{{\SetFigFont{6}{7.2}{\rmdefault}{\mddefault}{\updefault}p}}}}}
\put(2712,7749){\makebox(0,0)[lb]{\smash{{{\SetFigFont{6}{7.2}{\rmdefault}{\mddefault}{\updefault}p'}}}}}
\put(5037,5274){\makebox(0,0)[lb]{\smash{{{\SetFigFont{6}{7.2}{\rmdefault}{\mddefault}{\updefault}r'}}}}}
\put(4137,6024){\makebox(0,0)[lb]{\smash{{{\SetFigFont{6}{7.2}{\rmdefault}{\mddefault}{\updefault}c(A)}}}}}
\put(2712,6549){\makebox(0,0)[lb]{\smash{{{\SetFigFont{6}{7.2}{\rmdefault}{\mddefault}{\updefault}c(A')}}}}}
\put(4962,4749){\makebox(0,0)[lb]{\smash{{{\SetFigFont{6}{7.2}{\rmdefault}{\mddefault}{\updefault}t'}}}}}
\put(4512,3699){\makebox(0,0)[lb]{\smash{{{\SetFigFont{6}{7.2}{\rmdefault}{\mddefault}{\updefault}c(B)}}}}}
\put(6087,3399){\makebox(0,0)[lb]{\smash{{{\SetFigFont{6}{7.2}{\rmdefault}{\mddefault}{\updefault}q}}}}}
\put(2637,2124){\makebox(0,0)[lb]{\smash{{{\SetFigFont{6}{7.2}{\rmdefault}{\mddefault}{\updefault}q'}}}}}
\put(2712,3324){\makebox(0,0)[lb]{\smash{{{\SetFigFont{6}{7.2}{\rmdefault}{\mddefault}{\updefault}c(B')}}}}}
\end{picture}
}

\end{center}
\medskip

But $A$- and $B$-domains intersect only in point $s$, because the lines
through $(c(A), p)$ and $(s, r')$ diverge on the hyperbolic plane. In fact, denote by
$ \alpha_1$, $ \alpha_2$, $ \beta_1$, $ \beta_2$ the angles
(p, c(A), s), (c(A), s, r'), (c(A), p, r'), (p, r', s), respectively. They
are equal to $\frac{4 \pi}{m} + \frac{2 \pi}{m}$, $\frac{\pi}{m}$,
$\frac{\pi}{m} + \frac{\pi}{m}$, $\frac{2 \pi}{m} + \frac{\pi}{m}$,
respectively. So $\alpha_1 + \alpha_2 = \frac{7 \pi}{m} \le \pi$, since 
$m \ge 7$ and the lines, if they converge or parallel, do it on the right
side of Fig. 1b. Now, $\beta_1 + \beta_2 = \frac{5 \pi}{m} < \pi$ and the
lines, if they converge or parallel, do it on the left. So, they diverge.

We demonstrated ad absurdum, the non-existence of the vertex $s*$
and so, of the { \em closed} path $R$. So, a path $R$ of length $2k - 2$
is not closed. But $p, q$ is never an edge; so we need at least two edges in
order to close $R$. If two edges are enough, then points $r'$, $t'$ coincide,
i.e. $i = 2$; actually, two edges will be enough in the case $m = 7$.
The proof of Lemma 3 is completed.

\noindent
{\bf Proof of Theorem 2}

Consider star-$m$-gons $A$, $B$ and the circuit $C$ as in beginning of the
 proof
of Lemma 3 above. Take the arc $e = (a_0, a_1)$ on the circuit $C$; by Lemma 1
(i), $e$ is balanced,
i.e. the vectors $x_i (e^*)$ of the opposite arc $e^* = (b_0, b_1)$ are the 
same, as of the arc $e$, but they have opposite directions with
respect of the circuit $C$. The same arc $e$, seen as an arc of the circuit $B$
of length $m$, is opposite to two arcs in this odd circuit and, in particular,
to the arc $(a_k,a_{k +1})$. The last arc has, by
Lemma 1 (ii), $ \frac{ \lambda}{2}$ vectors, coinciding
with vectors of $e$, but with opposite direction on the circuit $B$. Finally,
consider the oriented path $(a_{k +1},
a_k = b_0, b_1)$ of length 2 in our $\{ \frac{m}{2}, m\}$. Its two
arcs have vectors, coinciding, but going in opposite direction on this path.
But it contradicts to Lemma 2, because $2 < k  $.

\section{ Spherical analogue of Coxeter's honeycombs}

In this Section we consider, for any pair $(i,m)$ of integers, such that
$1 \le i < \frac{m}{2}$ and $g.c.d.(i,m) = 1$, star-polygons $ \frac{m}{i}$.
Clearly, $ \frac{m}{1}$ denotes now a convex $m$-gon; so we see star-polygons
as a generalization of convex ones. We will allow further extension:
star-polygons $ \frac{m}{i}$ with $ \frac{m}{2} < i < m$, let us call them
{\em large}  star-polygons. They cannot be represented on Euclidean or
hyperbolic plane, because they have there the same representation as 
$ \frac{m}{m-i}$.
But they can be represented on the sphere by the following way; see Fig. 2 for
the simplest $\frac{3}{1}$ and $\frac{3}{2}$. Let
$a_0,...,a_{m-1}$ be $m$ points, placed in this order, on a great circle of the
sphere, in order to form a regular $m$-gon. Then the spherical (great circle)
distance $d(a_0, a_i)$ is $ \frac{2 \pi i}{m}$, but on $ \frac{m}{i}$, the
length of the way is $d(a_0,a_i)$ for $i < \frac{m}{2}$ and
$2 \pi - d(a_0,a_i)$ otherwise. Using this larger set of polygons, we will
look for spherical representations of regular (i.e. with a group of symmetry
acting transitively on all $j$-faces, $0 \le j \le 2$) polyhedra.

\medskip
\begin{center}
\setlength{\unitlength}{0.00058333in}
\begingroup\makeatletter\ifx\SetFigFont\undefined%
\gdef\SetFigFont#1#2#3#4#5{%
  \reset@font\fontsize{#1}{#2pt}%
  \fontfamily{#3}\fontseries{#4}\fontshape{#5}%
  \selectfont}%
\fi\endgroup%
{\renewcommand{\dashlinestretch}{30}
\begin{picture}(9016,4245)(0,-10)
\put(5452.771,39.334){\arc{11095.774}{3.3347}{3.9956}}
\put(-1836.772,39.334){\arc{11095.775}{5.4292}{6.0901}}
\put(7170.500,831.272){\arc{3696.216}{3.2488}{6.1760}}
\put(1808.000,6354.000){\arc{11100.000}{1.2405}{1.9011}}
\put(6080.354,2719.616){\arc{3696.875}{5.3435}{8.2702}}
\put(1208,54){\makebox(0,0)[lb]{\smash{{{\SetFigFont{8}{9.6}{\rmdefault}{\mddefault}{\updefault}   Fig.  2a.     3/1}}}}}
\put(8260.902,2719.295){\arc{3696.080}{1.1546}{4.0810}}
\put(6608,54){\makebox(0,0)[lb]{\smash{{{\SetFigFont{8}{9.6}{\rmdefault}{\mddefault}{\updefault}    Fig.  2b.     3/2}}}}}
\end{picture}
}

\end{center}
\medskip

In the Table 2 below, the rows (columns) denote a cell (respectively, a
vertex figure) of would-be representations. If the representation,
corresponding
to a given pair of $( \frac{m}{i}, \frac{n}{j})$ of polygons, exists, we
denote it by this pair and write its density in corresponding cell of the
Table 2. The densities were counted directly, by superposing the representation
on corresponding regular polyhedron. But the expression of the density, given
in the formula 6.41 of \cite{Cox1} for multiply-covered sphere is valid
for our representations, i.e. the density of $( \frac{m}{i}, \frac{n}{j})$
is $N_1 ( \frac{i}{m} + \frac{j}{n} - \frac{1}{2})$, where $N_1$ is the number
of edges. (Above expression is equivalent to Cayley's generalization of 
Euler's Formula, given as the formula 6.42 in \cite{Cox1}.)
Our representations are {\em Riemann surfaces}, i.e. $d$-sheeted spheres (or
 $d$ almost coincident, almost spherical
surfaces) with the sheets connected in certain branch-points.  

We see a $ \frac{m}{i}$ as a representation of the $m$-cycle on the sphere,
together with a bi-partition of $i$-covering of the sphere. Call
{ \em interior} the part with angle, which is less than $\pi$. For
representations
below, the vertex figure selects uniquely the part of the cell: namely,
the vertex figure $ \frac{n}{j}$ gives the value $ \frac{2 \pi j}{n}$ for
the angle of the cell. It takes interior of the cell if $j < \frac{n}{2}$
and exterior otherwise.
 
The Table 2 shows that each of all nine regular polyhedra (seen as abstract
surfaces) admits four
such Riemann surfaces and we checked, case by case, that all 36 are different
and that remaining 28 possible representations do not exist. Each  of four
representations for each regular polyhedron has same genus as corresponding
abstract surface; so the genus is four for 8 representations of the form
$( \frac{5}{i}, \frac{5}{j})$ and zero for all others. 

In the Table 2, the column with $ \frac{2}{1}$ corresponds to
doubling of regular polygons. Alexandrov ( \cite{Al}) considered, for
other purpose, the doubling of any
convex polygon as an abstract sphere, realized as a degenerated (i.e. with
volume 0) convex polyhedron. $m2$ and $2n$ on the plane and the sphere
appeared also in Section 7 of \cite{FT}.  By analogy, we will do such doubling
for star-polygons $ \frac{m}{i}$ with $i < \frac{m}{2} $. But for large 
star-polygons we should do doubling on the sphere. The row and the column
with $\frac{m}{i}$ correspond to any pair of mutually prime integers
$(i,m)$, $1 \le i < m $. As Table 2 shows, there exist all representations 
$( \frac{2}{1}, \frac{m}{i})$ and $( \frac{m}{i}, \frac{2}{1})$ and each of
them has density $i$ (and the genus 0).

An infinity of other representations can be obtained by permitting polygons
$ \frac {m}{i + tm}$ for any integer $t \ge 0$; the way on the edge 
$(a_0, a_{i + tm})$ will be $ 2 \pi t - d(a_0, a_{i + tm})$.

{ \bf Table 2.
 36 representations of regular polyhedra on the sphere. }
\[\begin{array}{|c||c||c|c||c|c||c|c||c|c||c|c||} \hline
& \frac{2}{1}& \frac{3}{1}& \frac{3}{2}& \frac{4}{1}& \frac{4}{3}& \frac{5}{1}& \frac{5}{4}& \frac{5}{2}& \frac{5}{3}& \frac{m}{i}& \frac{m}{m-i} \\ \hline \hline
 \frac{2}{1}&1&1&2&1&3&1&4&2&3&i&m-i \\ \hline \hline
 \frac{3}{1}&1&1&3&1&7&1&19&7&13&& \\ \hline
 \frac{3}{2}&2&3&5&5&11&11&29&17&23&& \\ \hline \hline
 \frac{4}{1}&1&1&5&&&&&&&&\\ \hline
 \frac{4}{3}&3&7&11&&&&&&&&\\ \hline \hline
 \frac{5}{1}&1&1&11&&&&&3&9&&\\ \hline
 \frac{5}{4}&4&19&29&&&&&21&27&& \\ \hline \hline
 \frac{5}{2}&2&7&17&&&3&21&&&& \\ \hline
 \frac{5}{3}&3&13&23&&&9&27&&&& \\ \hline \hline
  \frac{m}{i}&i&&&&&&&&&& \\ \hline
  \frac{m}{m-i}&m-i&&&&&&&&&& \\ \hline \hline
\end{array} \]

\section{ Star-honeycombs}

Besides star-polygons and four regular star-polyhedra on 2-sphere, which are
all embeddable (last four are isomorphic to Ico or Do), there are
( \cite{Cox}) only following regular star-honeycombs: ten regular
star-polytopes on 3-sphere and four star-honeycombs
in hyperbolic 4-space; see the Tables 1, 3-5. In this Section we show that none
of last 14 is embeddable. Consider first the case of 3-sphere.

There are six regular 4-polytopes (4-simplex $ \alpha_4$, 4-cross-polytope
$ \beta_4$, 4-cube $ \gamma_4$, self-dual 24-cell and the pair of dual
600-cell and 120-cell) and ten star-4-polytopes; see the Chapter 14 in 
 \cite{Cox1}.
 \cite{As} showed non-embeddability of 24- and 600-cell; \cite{DG3} did it
for 120-cell. Clearly, $ \gamma_4$ and $ \beta_4$ are $H_4$ and
$ \frac{1}{2} H_4$ themselves.

Embeddable ones among Archimedean tilings of 3-sphere and 3-space, were
identified in  \cite{DS4}; for example, { \em snub} 24-cell (semi-regular
Gosset's 4-polytope $s\{3, 4, 3\}$) embeds into $ \frac{1}{2} H_{12}$ while
the Grand Antiprism of  \cite{Con} is not embeddable.

The isomorphisms among ten star-4-polytopes, see \cite{vO} and  pages 266-267
of \cite{Cox1}, preserve all incidencies and imply, of course, isomorphisms
of the skeletons of those
polytopes. Using Schl\"afli notation, those isomorphisms {\em of graphs} are:

(i) $\{\frac{5}{2},5,3\} \sim \{5, \frac{5}{2}, 3\}$;

(ii) $\{ \frac{5}{2},3.3\} \sim $120-cell (remind the isomorphism of
$ \{\frac{5}{2}, 3\}$ and $\{5,3\}$);

(iii) all remaining seven skeletons are isomorphic with the skeleton of
 600-cell (moreover, $\{3,5, \frac{5}{2}\}$ has same faces; remind the isomorphism
of $\{3, \frac{5}{2}\}$ and $\{3,5\}$).

So eight star-polytopes from (ii) and (iii) above are not embeddable. Remaining
case (i) is decided by the Theorem 3 below, using 5-gonal inequality.

\begin{theor}
\label{4star}
None of ten star-4-polytopes is embeddable.

\end{theor}

\noindent
{\bf Proof of Theorem 3}\

In view of above isomorphisms, it will be enough to show that (the skeleton
of) 4-polytope $P:= \{\frac{5}{2},5,3\}$ is not 5-gonal. P is the stellated
120-cell and $ \{\frac{5}{2}, 5\}$ is the (small) stellated dodecahedron, i.e. all
face-planes are extended until their intersections form a pyramid on each
face. P has 120 vertices, as 600-cell; namely, the centers of all 120 
(dodecahedral) cells of 120-cell. For any vertex $s$ of $P$, denote by Do(s)
the corresponding dodecahedron. $P$ has (as 120-cell) 1200 edges, 720
faces and 120 cells; its density is 4. Any edge $(s,t)$ of $P$ goes through
interiors of $Do(s)$, $Do(t)$ and the edge of 120-cell, linking those
 dodecahedra; $(s,t)$ is a continuation of this edge in both directions till
 the centers of dodecahedra $Do(s)$,$Do(t)$.

Consider now Fig. 3. Take as vertices $a$ and $b$ (for future
contre-example for 5-gonal
inequality) some two vertices of $\{ \frac{5}{2}, 5\}$
(a cell of $P$), which are centers of two face-adjacent dodecahedral cells of
120-cell. Let $Q := (q_1, q_2, q_3,q_4, q_5)$ be this common face of adjacency,
presented by the 5-cycle of its vertices. For any $q_i$ there is unique
star-5-gon $(a, d_i, b, d'_i, d''_i)$, such that sides $(b, d'_i)$ and
$(d''_i, a)$ intersect in the point $q_i$. Now, $D:=(d_1, d_2, d_3, d_4,
d_5)$ is a 5-cycle in $P$, because each $(d_{i-1},d_i)$
is an edge  in one of five cells $\{ \frac{5}{2}, 5\}$ of $P$, containing vertices
$a$ and $b$.  Put $x:=d_1$, $y:=d_2$, $z:=d_4$ and
check that the 5-gonal inequality for five vertices $a, b, x, y, z$ of $P$,
does not hold.

\medskip

\begin{center}
\hspace*{35mm}\epsfig{file=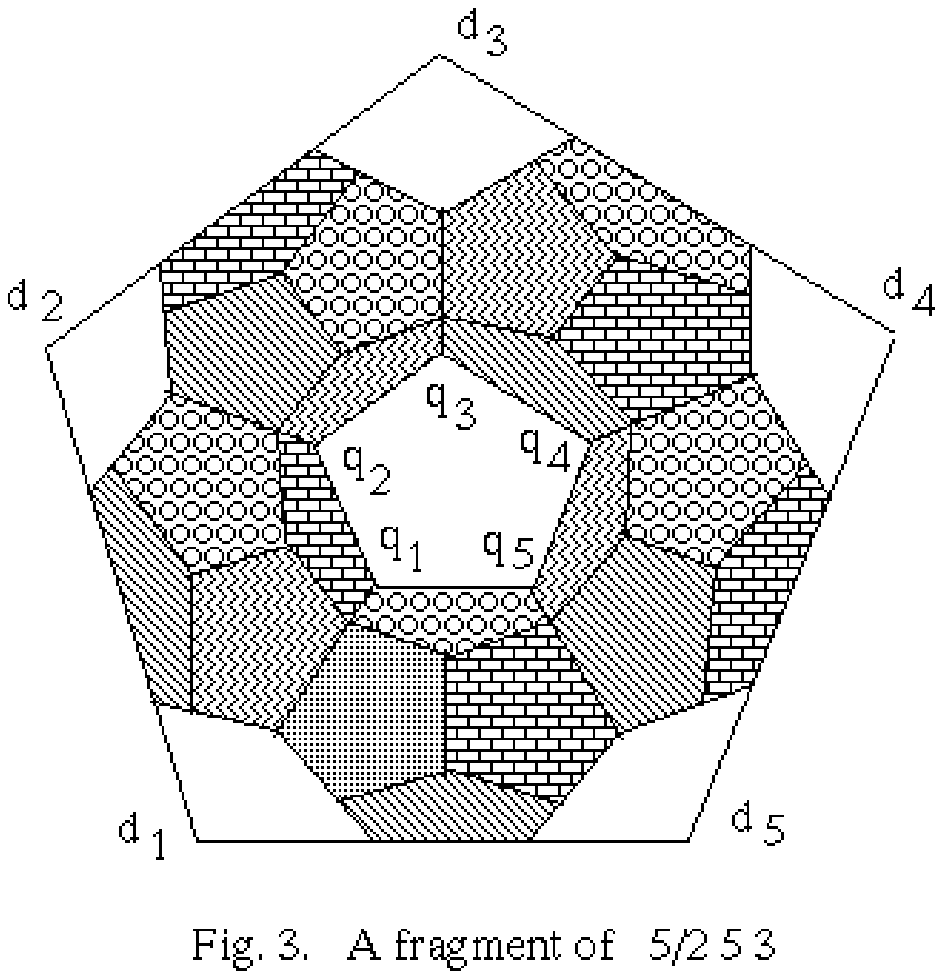,height=120mm} 
\end{center}

\medskip

In fact, $d_{xy} = 1 = d_{ax} = d_{ay} = d_{az} = d_{bx} = d_{by} = d_{bz}$,
because of the presence of corresponding edges in $P$. Therefore,
$d_{xz}, d_{yz}$ and $d_{ab}$ are at most 2. So, the absence of edges (x,z),
(y,z) and (a,b) will complete the proof of the Theorem 3. The edge $(a,b)$ does
not exist, because $Do(a)$ is face-adjacent to $Do(b)$. The edge
$(x,z)$ does not exists, because the line, linking vertices $x$ and $z$, goes,
 besides $Do(x)$
and $Do(z)$, through two other dodecahedra (such that their stellations are
$\{ \frac{5}{2}, 5\}$, containing vertices $a,b,d_2,d_3$ or $a,b,d_3,d_4$).
By symmetry, the edge $(y,z)$ does not exist also. We are done.

\noindent
{\bf Corollary}\ {\em None of four star-honeycombs in hyperbolic 4-space is
embeddable.}

\noindent
{\bf Proof of Corollary}\

In fact, $ \{\frac{5}{2},5,3,3\}$ has cell which contains (because of the
Theorem 3), as an induced subgraph, non-5-gonal graph $K_5 - K_3$. But any
induced graph of diameter 2 is isometric; so $\{ \frac{5}{2},5,3,3\}$ is not
5-gonal. $\{3,3,5, \frac{5}{2}\}$ has cell  $\{3,3,5\} =$ 600-cell. Two other hiiave cells which are
isomorphic to 600-cell. But 600-cell (seen by Gosset's construction as 
capping of all 24 icosahedral cells of snub 24-cell) contains also a forbidden
 induced graph
of diameter 2: pyramid on icosahedron (it violates {\em 7-gonal} inequality,
which is also necessary for embedding; see \cite{D}, \cite{DS1}). So, three
other star-4-polytopes are also non-7-gonal and non-embeddable.

\section{Regular tilings of dimension $d \ge 3$ }

The Tables 3-5 below present all of them and also all regular honeycombs
in the dimensions 3, 4, 5; for higher dimensions, $(d+1)$-simplices
$ \alpha_{d+1}$, $(d+1)$-cross-polytopes $ \beta_{d+1}$, $(d+1)$-cubes
$ \gamma_{d+1}$ and cubic lattices $ \delta_d$ are only regular ones.

In those Tables, $24-$, $600-$, $120-$ are regular spherical 4-polytopes
$\{3,4,3\}$, $\{3,3,5\}$, $\{5,3,3\}$ with indicated number of cells and $De(D_4)$, $Vo(D_4)$
are regular partitions $\{3,3,4,3\}$, $\{3,4,3,3\}$ of Euclidean 4-space, which are also
Delaunay (Voronoi, respectively) partitions associated with the (point) lattice
$D_4$.

All cases of embeddability are marked be the star * in the Tables. As in
Table 1 above, we omit in Tables 3-5 (in order to fit them in the page) the
brackets and commas in Schl\"afli notation.

{\bf Table 3. 
 3-dimensional regular tilings and honeycombs.}
\[\begin{array}{|c||c|c|c|c|c||c|c|c||c|c|c|c||} \hline
 & \alpha_3 &\gamma_3 & \beta_3 & Do & Ico& \delta_2& 
63&36&3\frac{5}{2}&\frac{5}{2}3&5 \frac{5}{2}&\frac{5}{2} 5\\ \hline \hline 
\alpha_3 & \alpha_4* & &\beta_4* & &600-&&&336 &33 \frac{5}{2}& && \\ \hline
 \beta_3 &&24- &&&&344 &&& &&&    \\ \hline
\gamma_3  &\gamma_4*  &&{\bf \delta_3*}  &  & 435*   & & &436*&&&& \\ \hline
Ico &  && & 353 & &&   && &  & 35\frac{5}{2}& \\ \hline
Do &120-    &   &534 &&535& & &  536&53 \frac{5}{2}& && \\ \hline \hline
\delta_2& &443* &&& &444*  && &  &&& \\ \hline
36 &&& & && &363& & & && \\ \hline
63 &633* & &634* &&635* & &&636* &&&& \\ \hline \hline
\frac{5}{2} 3&\frac{5}{2} 33 & & &&\frac{5}{2} 35&&&&&&& \\ \hline
3 \frac{5}{2}& & & && & &  &&&&&3 \frac{5}{2} 5 \\ \hline 
\frac{5}{2} 5&&&&\frac{5}{2} 53&&&&&&&\frac{5}{2} 5 \frac{5}{2} & \\
\hline
5 \frac{5}{2}&&&&&&&&&&5 \frac{5}{2} 3&&5 \frac{5}{2} 5 \\ \hline\hline
\end{array} \]

{\bf Table 4. 
 4-dimensional regular tilings and honeycombs.}
\[\begin{array}{|c||c|c|c|c|c|c||c||c|c|c||} \hline
 & \alpha_4 &\gamma_4 & \beta_4 & 24- & 120-& 600-&\delta_3& 
35\frac{5}{2}&\frac{5}{2}53&5\frac{5}{2}5\\ \hline \hline 
\alpha_4 & \alpha_5* & &\beta_5* & && 3335 && && \\ \hline        
\beta_4 && &&{\bf De(D_4)} &&& &&&    \\ \hline
\gamma_4  &\gamma_5*  && {\bf \delta_4*}  &  & & 4335*   & & && \\ \hline
24- &   & {\bf Vo(D_4)} &    && &  & 3434&&& \\ \cline{3-5}\hline
600- &    &   & &&& & &  335\frac{5}{2}& & \\ \hline
120- &5333 & &5334  && & 5335&  &&& \\ \hline \hline
\delta_3 & & & &4343*& & & &&& \\ \hline \hline
\frac{5}{2}53 & & & &&\frac{5}{2}533 & && && \\ \hline
35\frac{5}{2}& & & && && &&&35\frac{5}{2}5 \\ \hline
5\frac{5}{2}5& & & && & &  &&5\frac{5}{2}53& \\ \hline \hline
\end{array} \]

{\bf Table 5.
 5-dimensional regular tilings and honeycombs.}
\[\begin{array}{|c||c|c|c||c|c|c||} \hline
 &\alpha_5&\gamma_5&\beta_5&Vo(D_4)&De(D_4)&\delta_4 \\ \hline \hline
\alpha_5&\alpha_6*&&\beta_6*&&& \\ \hline
\beta_5&&&&&33343& \\ \hline
\gamma_5&\gamma_6*&&{\bf \delta_5*}&&& \\ \hline \hline
De(D_4)&&&&33433&& \\ \hline
Vo(D_4)&&34333&&&&34334 \\ \hline
\delta_4&&&&&43343*& \\ \hline \hline
\end{array}{} \]

Theorems 1, 2 above show that all regular 2-dimensional tilings and
star-honeycombs are embeddable except $\{ \frac{m}{2},m\}$ for all odd $m \ge 7$.
The following Theorem decides all remaining regular cases.

\begin{theor}
\label{dim4}

All embeddable regular tilings and honeycombs of dimension $d \ge 3$ are
 tilings:

(i)  either $ \alpha_{d+1}$, or $ \beta_{d+1}$, or

(ii) all with bipartite skeleton:

  (ii-1) all with cell $ \gamma_d$: $ \gamma_{d+1}$, $ \delta_d $ and
3 hyperbolic ones: $\{4,3,5\}$, $\{4,3,3,5\}$, non-compact $\{4,3,6\}$;

  (ii-2) all 4 with cell $ \delta_{d-1}$: hyperbolic non-compact $\{4,4,3\}$,
$\{4,4,4\}$, $\{4,3,4,3\}$ and $\{4,3,3,4,3\}$;

  (ii-3) all 4 with cell $\{6,3\}$: hyperbolic non-compact $\{6,3,3\}$,
$\{6,3,4\}$, $\{6,3,5\}$, $\{6,3,6\}$.

All $l_1$-rigid regular tilings are the bipartite ones; all
bipartite ones (except $ \gamma_{d+1}$ and $ \delta_d$ themselves ) embed
into $ { \bf Z}_ \infty $. 

\end{theor}

\noindent
{\bf Proof of Theorem 4}

In fact, we review all cases of Tables 3-5. All compact cases (on first 5
rows, columns of Table 3 and first 6 rows, columns of Table 4) were decided
in \cite{DS2}. Non-embeddability for all 14 star-polytopes and
star-honeycombs (in Tables 3, 4) was established in Section 4.
It remains 11, 2, 5 non-compact tilings of hyperbolic 3-, 4-, 5-space; we will
show that 7, 1, 1, respectively, of them are embeddable into
${ \bf Z}_ \infty$, while 8 others are not 5-gonal.  

The tilings $\{3,4,3,4\}$, $\{3,4,3,3,3\}$, $\{3,3,4,3,3\}$, $\{3,4,3,3,4\}$
have non-5-gonal graph
$K_5 - K_3$ as induced subgraph of the cell. $\{3,6,3\}$ (respectively,
$\{3,4,4\}$)
contain induced $K_5 - K_3$, because each its edge is common to 3
(respectively, to 4) triangles. $\{3,3,6\}$ is a simplicial manifold with 6
triangles on an edge; taking 1-st, 3-rd and 5-th of them, we get again induced
$K_5 - K_3$. A particularity of $T:=\{3,3,3,4,3\}$ is that the cell $ \beta_4$ of its
vertex figure $De(D_4)$ is also the equatorial section of 
the cell $ \beta_5$ of $T$. All neighbors of a vertex $s$ of $T$ form
 $De(D_4)$. Take an isometric subgraph $K_5 - K_3 $ in $De(D_4)$, given in
 \cite{DS3}. The vertex $s$ is a neighbor of each of its five vertices;
obtained 6-vertex graph is non-5-gonal graph of diameter 2, which is, using
above particularity of $T$, is an induced subgraph of $T$. (Compare with 
embeddable tiling $\{4,3,3,4,3\}$ by $\gamma_5$, having the same vertex figure.) All
seven above
tilings are not 5-gonal, because any induced graph of diameter 2 is isometric.
Finally, each edge of $\{5,3,6\}$ is common to 6 disjoint pentagons; taking
1-st, 3-rd and 5-th of them we obtain non-5-gonal 11-vertex induced subgraph
of diameter 4 of $\{5,3,6\}$; a routine check shows that it is isometric.  

Other hyperbolic tilings embed into ${ \bf Z}_ \infty$, because of Lemma 5
below; it is easy to find reflections, required by Lemma 5 in each case.
It is easy to check $l_1$-rigidity for all (except of Tetrahedron, which is
not $l_1$-rigid) cases of embedding for dimension 2. Now,
any bipartite  embeddable graph is $l_1$-rigid, because it has scale 1. The
proof is complete.

Let $T$ be any (not necessary regular) convex $d$-polytope or tiling of
Euclidean or hyperbolic $d$-space by convex polytopes, such that the
skeleton is a bipartite graph. (We admit infinite cells and, if regular,
infinite vertex figures.) Then the set of its edges can be partitioned
into {\em zones}, i.e.
sequences of edges, such that any edge of a sequence is the opposite to the
previous one on a 2-face (which should, therefore, be even). 

\noindent

{\bf Lemma 5}
{\em Let $T$ is as above; suppose that the mid-points of edges
of each zone lie on hyperplanes, different for each zone, which are (some of)
reflection hyperplanes of $T$ and perpendicular to edges of their zones. Then
$T$ embeds into ${ \bf Z}_m$ with $m$ no more than the number of zones. }

\noindent

{\bf Proof of Lemma 5}

It follows directly from the fact that each geodesic path (in the skeleton of
$T$) intersects any zone in at most one edge.

{\bf Remark 3}

Embedding of any bipartite regular tiling can be obtained, using Lemma 5.
The reflections, required by Lemma 5 (let us call them {\em zonal} reflections)
generate, because of simple connectedness of $T$, a vertex-transitive group
of automorphisms of $T$ (call it {\em zonal} group); so $T$ is uniform and the
zonal group is generated by the zonal reflections of all edges incident to
a fixed vertex of $T$.
For any  fixed 2-face of $T$, which is a $2k$-gon, let $m_1$,...,$m_k$  be the
zonal reflections of its edges, considered in the cyclic order. Then the
product $m_{1}...m_{k}m_{1}...m_{k} = <1>$ (i.e. $m_1...m_k$ is an
involution) and those relations, for all 2-faces around a fixed vertex of $T$,
are all defining relations for the zonal group of $T$. So, the zonal group is 
not 2-transitive on vertices. For example, the zonal group of
Archimedean truncated $\beta_3$
is an 1-transitive subgroup of index 2 of the octahedral group
$Aut(T)=O_h$, which is 2-transitive. Also, a polytope in the
conditions of Lemma 5 is not necessary zonotope. For example, any
centrally-symmetric non-Archimedean (by choice of the length of truncation)
truncated $\beta_3$ fits in it; it is a zonohedron in original sense of
Fedorov, but not in usual sense of Minkowski (with all edges of each zone 
having same length).

{\bf Remark 4}

 All infinite families of regular tilings are embeddable. In fact, $m$-gons,
$ \delta_{n-1} = { \bf Z}_n$, $ \gamma_n = H_n $, $ \alpha_n$, $ \beta_n$ are
embeddable and, moreover, first three are $l_1$-rigid.
But embeddings of skeletons of $\alpha_n$ and, for $n \ge 4$, $\beta_n$, is
more complicate. It is considered in detail (in terms of corresponding
complete graph $K_{n+1}$ and Cocktail-Party graph $K_{n \times 2}$ in 
Chapter 23 \cite{DL} and Section 4 of Chapter 7 \cite{DL}, respectively.
Any $ \alpha_n, n \ge 3$, is not $l_1$-rigid, i.e. it admits at least two
different embeddings. We give now two embeddings of $ \alpha_n $ into
 $m$-cubes with scale  $ \lambda $ , realizing, respectively, maximum and
minimum of $ \frac {m}{ \lambda } $. The first one is 
 $ \alpha_n \to \frac {1}{2}
H_{n+1}$. Now define  $m_n = \frac {2n}{n+1}$  for odd  $n$  and
 $ = \frac {2n+2}{n+2} $  for even  $n$; define $ \lambda_n $ be the minimal
even positive number  $t$  such that  $tm_n$  is an integer. Then 
 $ \alpha_n $ embeds into $tm_n$-cube with scale  $ \lambda_n $; for example,
 $ \alpha_4$ embeds into 10-cube with scale 6. Any
 $ \beta_n $, $n \ge 4$, is not $l_1$-rigid. All embeddings of  $ \beta_n $
are into $ 2 \lambda $-cube with any such even scale $ \lambda $ that
$ \alpha_{n-1} $ embeds into $m$-cube, $ m \le 2 \lambda $, with scale 
$ \lambda $. For minimal such scale, denote it  $ \mu_n $, the following
is known:  $ n > \mu_n \ge 2 \lceil \frac {n}{4} \rceil $  with equality
 in the lower bound for any  $n \le 80 $ and, in the case of $n$
divisible by 4, if and only if there exists an Hadamard matrix of order  $n$.
In particular,  $ \beta_3 \to \frac {1}{2} H_4 $, $ \beta_4 \to \frac {1}{2}
H_4 $ (in fact, they coincide as 4-polytopes, but there are two embeddings),
and $ \beta_5 $ embeds only with scale 4 (into  $ H_8$).

{\bf Remark 5}

This note finalizes the study of embeddability for regular tilings done 
in \cite{DS1}, \cite{DS2}; we correct now following misprints there:
a) in the sentence ``Any $l_1$-graph, not containing $K_n$, is
$l_1$-rigid'' on p.1193 \cite{DS1}, should be $K_4$ instead of $K_n$;
b) in the sentence, on p.1194 \cite{DS1}, about partitions of Euclidean plane,
embeddable into ${\bf Z}_m$, $m < \infty$, should be $ \le $ instead of $ < $; 
c) in the sentence about F\"oppl partition on p.1292 \cite{DS2}, should be
$\alpha_3$ and truncated $\alpha_3$ instead of $\alpha_3$.

\end{document}